\documentclass[11pt]{amsart}

\usepackage{algorithm}
\usepackage[noend]{algpseudocode}
\usepackage{amsfonts}
\usepackage{amsmath}
\usepackage{amssymb}
\usepackage{amsthm}
\usepackage[toc,page]{appendix} 
\usepackage{array}
\usepackage{bbm}
\usepackage{colortbl} 
\usepackage{enumerate}
\usepackage{eurosym}
\usepackage{float}
\usepackage[tmargin=2.0cm, bmargin=2.0cm, lmargin=2.0cm, rmargin=2.0cm]{geometry}
\usepackage{graphicx}
\usepackage{graphics}
\usepackage{latexsym}
\usepackage{longtable}
\usepackage[section]{placeins}
\usepackage{setspace}
\usepackage{sidecap}
\usepackage{siunitx}  
\usepackage{subfig}
\usepackage{wrapfig}
\usepackage[table]{xcolor} 
\usepackage{amsaddr}

\newcolumntype{L}[1]{>{\raggedright\let\newline\\\arraybackslash\hspace{0pt}}m{#1}}
\newcolumntype{C}[1]{>{\centering\let\newline\\\arraybackslash\hspace{0pt}}m{#1}}
\renewcommand{\arraystretch}{1.2}

\makeatletter
\def\maxwidth{ %
  \ifdim\Gin@nat@width>\linewidth
    \linewidth
  \else

    \Gin@nat@width
  \fi
}
\makeatother

\definecolor{fgcolor}{rgb}{0.345, 0.345, 0.345}

\usepackage{framed}
 {\par\unskip\endMakeFramed%
 \at@end@of@kframe}
\makeatother

\definecolor{shadecolor}{rgb}{.97, .97, .97}
\definecolor{messagecolor}{rgb}{0, 0, 0}
\definecolor{warningcolor}{rgb}{1, 0, 1}
\definecolor{errorcolor}{rgb}{1, 0, 0}
\newenvironment{knitrout}{}{} 

\title{An application of dynamic programming \\to assign pressing tanks \\at wineries}

\author{Zbigniew Palmowski}
\address{Faculty of Pure and Applied Mathematics, Wroc\l aw University of Science and Technology,\\ ul. Wyb. Wyspia\'nskiego 27, 50-370 Wroc\l aw, Poland}
\email{zbigniew.palmowski@gmail.com}

\author{Aleksandra Sidorowicz}
\address{Faculty of Pure and Applied Mathematics, Wroc\l aw University of Science and Technology,\\ ul. Wyb. Wyspia\'nskiego 27, 50-370 Wroc\l aw, Poland}
\email{sidorowicz.aleksandra@gmail.com}

\thanks{Corresponding author: Zbigniew Palmowski.}
\thanks{This work is partially supported by the National Science Centre, Grant No. 2016/23/B/HS4/00566 (2017-2020).}

\date{\today}

\begin{document}

\begin{abstract}
This paper describes an application of dynamic programming  to determine the optimal strategy for assigning grapes to pressing tanks in one of the largest Portuguese wineries. To date, linear programming has been employed to generate proposed solutions to analogous problems, but this approach lacks robustness and may, in fact, result in severe losses in cases of sudden changes, which frequently occur in weather-dependent wine factories. Hence, we endowed our model with stochasticity, thereby rendering it less vulnerable to such changes. Our analysis, which is based on real-world data, demonstrates that the proposed algorithm is highly efficient and, after calibration, can be used to support winery's decision-making. The solution proposed herein could also be applied in numerous other contexts where production processes rely on outside supplies.

\vspace{3mm}

\end{abstract}

\maketitle

\pagestyle{myheadings} \markboth{\sc Z.\ Palmowski
--- A.\ Sidorowicz} {\sc Pressing tank assignment optimization}

\vspace{1.8cm}


\newpage

\section{Introduction}\label{chapt_1}

In this paper, we present an application of dynamic programming to determine the optimal scheduling policy for a grape-reception process, where the optimization criterion is profit maximization over a fixed time horizon. Similar problems affect many companies whose production processes rely on outside supplies and are very common, as supply chains are increasing in size and outsourcing is still very popular among general business strategies (see, for example, \cite{outsourcing}).

As \cite{12,13} point out, optimization procedures have taken on an increasingly important role in agricultural planning. Moreover, optimization is widely used in machine scheduling (see, e.g., \cite{8,2,1}), aircraft scheduling (see \cite{4,3}), vehicle assignments in delivery scheduling (see \cite{11}), gas purchase and transportation scheduling (see \cite{13}), and textile manufacturing (see \cite{14}). Lastly, \cite{5,6,7,15,9,10} adopt a more general perspective.

The specific optimization problem encountered in the context of a winery is more complex than those mentioned above, as deliveries of raw materials (i.e., grapes) can differ significantly from one another and, moreover, the winery can control neither the size of the delivery nor its timing. In particular, irregularities in time of delivery increase the intricacy of the problem and so necessitate the creation of an entirely new data-analysis methodology. The dynamic programming approach we have developed allows the decision-making process to adapt to the current state, which can vary unpredictably with respect to size, timing, and quality of grape harvest, and hence is better suited to the problem being addressed.

Our analysis involves the grape-reception process employed by one of the largest Portuguese companies producing Vinho Verde wine. Prior to reaching store shelves for purchase by consumers, Vinho Verde undergoes a complex process beginning with the grapes from which the wine is made to the finished wine’s bottling. The first stage is harvesting of the grapes, and their harvest date is a function of the crop's maturity, which is determined either by some empirical method or by laboratory tests based on concentration of sugars and total acidity. When the grapes are mature enough, harvesting can commence.

Two possible methods of picking grapes are manually, in which a worker employs shears to pick the fruit, and mechanically, wherein a mechanized grape harvester does so. The picked grapes are then transported from the vineyard to the winery, subject to certain timing and transport constraints selected to maximize the quality of the finished product, i.e., the wine. Firstly, pressing must occur as soon as possible after harvesting, as the grapes’ sweetness, acidity, and alcohol levels change rapidly following their harvest. Secondly, the grapes should be prevented from deteriorating by avoiding exposure to sunlight over an extended length of time and direct rainfall. Upon entering the destination winery, the truck conveying the grapes must be registered in a computerized system together with its arrival time, the variety and amount of grapes it is transporting. The truck then joins a queue as it awaits its turn for unloading  into one of the presses. When its turn to be unloaded comes, the truck’s load is emptied into a silo, where all pests are removed. This phase of the process ends with the grapes being cleaned; then they proceed to the pressing stage of the process, which the following steps comprise:  maceration (during which a freshly pressed fruit juice known as ‘must’ is produced); alcoholic fermentation of the must; malolactic fermentation of the must (in which bacteria convert tart-tasting malic acid, naturally present in grape must, into softer-tasting lactic acid); laboratory analysis; and finally bottling of the ready-to-drink wine.

A wine-producing company's most challenging optimization problem involves distribution of grape deliveries among the winery’s pressing tanks so as to obtain the highest possible profit  from the wine being produced. The primary difficulty posed by this assignment problem is the inability of the winery to control the delivery times of grape harvests, which typically originate with outside suppliers. In addition, there are four different types of grapes, which must be processed separately and which produce distinct types of wine that are sold for differing prices.

Furthermore, as we will demonstrate, the press assignment task can be modeled as a Markov control process. 
That in turn allows to apply the Bellman, or dynamic programming, equation, in which the problem is divided into more easily solved, overlapping optimization sub-problems, allowing identification of an optimal strategy (see, for example, Cormen et al. \cite[p. 359]{cormen}). To our knowledge, this specific real-time grape-to-press assignment problem has not been described extensively in the literature heretofore; instead, linear programming has been the optimization method most commonly used to solve grape-oriented  scheduling problems (see \cite{wine_supply_chain_overview,OR_methods_overview_fresh_fruits}), i.e., harvesting \cite{harvest_schedule_rising_market,robust_grapes,linear_grapes} and fermentation-tank management \cite{linear_tanks1,linear_tanks2,dec_support}. As can easily be seen, in our non-deterministic case, a linear approach is significantly less robust than a dynamic programming, which can be employed to solve the equivalent, above mentioned fermentation tank scheduling problems as well. Below we provide an extensive numerical analysis to illustrate the efficiency of the proposed algorithm.

The remainder of this paper is organized as follows. In Section \ref{chapt_2}, we introduce the Markov decision process and the Bellman equation to model the problem we seek to address. In Section \ref{sect_model_def}, we describe adapting a Markov decision process to assignment of grapes to pressing tanks. Section \ref{chapt_4} describes model's calibration employing the data available to us, and Section \ref{chapt_5} algorithm implementation. In Section \ref{chapt_6}, we present the results of applying the algorithm to numerous test cases. Finally, Section \ref{chapt_7} concludes the paper.

\section{Markov decision process}\label{chapt_2}
We first defined the temporal aspect of the Markov decision process; time $t$ within the process is discrete rather than continuous and has a finite rather than an infinite time horizon, i.e., time $t\in\{0,\ldots,T\}$, where $T$ is some finite positive number. The state of the stochastic system at time $t$ is $X_t$ and is defined by the following function:
\begin{equation*}\label{maindef}
X_t=f(t,X_{t-1},Y_t,Z_t).
\end{equation*}
Thus, $X_t$ is a function of the process’s state at time $t-1$, $X_{t-1}$; control $Y_t$ drawn from a specific set of allowed controls defined by $\Gamma(t,X_{t-1})$; and new information $Z_t$ arriving at time $t$. $Z_t$ is assumed to be a random variable independent of $Z_1,Z_2,\ldots,Z_{t-1}$, with a known probability distribution; that is, based on past information, the decision maker selects control $Y_t$, after which a piece of new and random information $Z_t$ arrives. The new state of the system at time $t$ is thus determined by function $f$ (known as the system’s transition function) containing the arguments presented above. During iterations at time $t\in \{1,\ldots,T\}$, the controller receives a payoff $g(t,X_{t-1},Y_t)$ generated from the previous state $X_{t-1}$ according to the current control $Y_t$.

The model’s goal is to obtain the highest possible sum of payoffs $g$ throughout the set of time intervals $\{1,\ldots,T\}$, which is approximated by the optimal value function $V^*$. The value function $V$ itself returns the expected sum of payoffs $g$ received by applying a certain strategy. Of course, the optimal value function $V^*$ is produced by an optimal strategy to select controls $Y_1^*, Y_2^*,\dots,Y_T^*$, which are denoted in this paper by $Y_{1\ldots T}^*$.
Note that the lower case symbols $x,y,z$ represent realizations of the equivalent upper case random variables.
Formally,
\begin{equation}\label{optimal_value_function}
V^*(t,x):=\sup_{Y_{t\ldots T}}V_{Y_{t\ldots T}}(t,x),
\end{equation}
where $Y_{t\ldots T}=(Y_t,\ldots,Y_T)$ and
\begin{equation*}
V_{Y_{t\ldots T}}(t,x):=\mathbb{E}\left[\sum_{s=t}^Tg(s,X_{s-1},Y_s)\middle|X_{t-1}=x\right].
\end{equation*}
Introducing notation $\mathbb{E}\left[\cdot \middle| X_{t-1}=x\right] := \mathbb{E}_x^{t-1}\left[\cdot\right]$,
\begin{equation}\label{value_function}
V_{Y_{t\ldots T}}(t,x):=\mathbb{E}_x^{t-1}\left[\sum_{s=t}^Tg(s,X_{s-1},Y_s)\right].
\end{equation}

It should be noted that the stochastic dynamical system $X$ can be modeled as a Markov process.
Hence, the state of the system at time $t$ should depend only on the prior state $X_{t-1}$, and not on any state occupied by the system prior to time $t-1$. By means of linearity of the conditional expected value and its tower property, it follows that
\begin{align}\label{transformation}
V_{Y_{t\ldots T}}(t,x)&=\mathbb{E}_x^{t-1}\left[g\left(t,X_{t-1},Y_t\right) \right]+\mathbb{E}_x^{t-1}\left[\sum_{s=t+1}^{T}g\left(s,X_{s-1},Y_s\right)\right] \nonumber\\
&=\mathbb{E}_x^{t-1}\left[g\left(t,X_{t-1},Y_t\right)\right]+\mathbb{E}_x^{t-1}\left[\mathbb{E}\left[\sum_{s=t+1}^{T}g\left(s,X_{s-1},Y_s\right)\middle| X_t\right]\right]\\
&=\mathbb{E}_x^{t-1}\left[g\left(t,X_{t-1},Y_t\right)\right]+\mathbb{E}_x^{t-1}\left[V_{Y_{t+1\ldots T}}\left(t+1,X_t\right)\right]\nonumber\\
&=\mathbb{E}_x^{t-1}\left[g\left(t,X_{t-1},Y_t\right)+V_{Y_{t+1\ldots T}}(t+1,X_t)\right]. \nonumber
\end{align}
Identity (\ref{transformation}) and Equations (\ref{optimal_value_function}) -- (\ref{value_function}) 
lead to the Bellman equation:
\begin{equation}\label{recursion_eq}
V^*(t,x)=\sup_{Y_t}\mathbb{E}_x^{t-1}\left[g\left(t,X_{t-1},Y_t\right)+V^*\left(t+1,X_t\right)\right],
\end{equation}
with the obvious boundary condition at the terminal time $T$:
\begin{equation}\label{final_value_eq}
V^*(T,x)=\sup_{Y_T}\mathbb{E}_x^{T-1}\left[g\left(T,X_{T-1},Y_T\right) \right].
\end{equation}

From the above, it is possible to compute the optimal value function $V^*(1,x_0)$ by means of backwards induction. The starting point is provided by Equation (\ref{final_value_eq}), which returns $V^*(T,x)$ for all possible $x$ values; then, using Equation (\ref{recursion_eq}) repeatedly, we can obtain $V^*(T-1,x),V^*(T-2,x),\ldots,V^*(2,x),V^*(1,x_0)$ sequentially. Important to note at this point is that only one component of a strategy is optimized at time $t$, not the entire strategy at once.

\section{Model adaptation}\label{sect_model_def}
We now describe modeling the winery's optimization problem as a Markov decision process. Thus, we define all model's components based on real-world data. First, recall that $T$ is the time horizon over which optimization is performed and that the Markov process is a discrete-time process.

{\bf State of stochastic dynamical system $\mathbf{X}$.}
Each $X_t = \{X_t^1,\ldots, X_t^K\}$, $t=0,\ldots,T$, represents the states of $K$ pressing tanks at time $t$, where $k=1,\ldots,K$, denotes a single pressing tank. In order to simplify the notation, we firstly assume that $K=1$, i.e., the manufacturer possesses only one pressing tank, and $X_t$ provides information regarding this single press. In this case, $X_{t}$ is described by a three-element vector $X_{t}=[v(X_t),l(X_t),s(X_t)]$, whose three elements are as follows: $v(X_t)$, the grape variety inside the one press; $l(X_t)$, weight of grapes;  and $s(X_t)$, the information on press  processing start time. Each of these is dependent on time $t$, what the notation specifying $t$ highlights. The process’s state $X_t$ depends on two other elements: press capacity $C$ and total processing time $TP$, both invariant over time. Note that $v(X_t)=l(X_t)=0$ indicates an empty press, and $s(X_t)=0$ is interpreted as a situation, when the press is not operating at time $t$ and so can be further filled.

{\bf Set of controls $\mathbf{Y}$.}
The main goal of this study is to identify  an optimal decision process. In fact, the controls $Y_t$ are interpreted as the decisions over which we are optimizing; namely, $Y_t$ assigns a truck  to the press being analyzed at time $t$; i.e., it picks the truck that should be unloaded at this point in time. The truck is selected from a queue that is filtered by the allowable set of truck types permitted to be unloaded at time $t$, $\Gamma(t,X_{t-1})$. Here,  “truck type” represents a combination of grape variety and weight contained in the truck's load, $(v,l)$; these will be defined in greater detail below. Each selected truck $Y_t=[v(Y_t),l(Y_t),a(Y_t)]$ is represented by a vector containing three elements: grape variety ($v(Y_t)$), load weight  ($l(Y_t)$), and arrival time ($a(Y_t)$), analogous to the notation employed with respect to state $X_t$ of the stochastic dynamical system.

As mentioned in Section \ref{chapt_1}, once grapes are harvested, their quality decreases with time, thus necessitating wineries to track grapes' times of arrival $a(Y_t)$ at the winery. When the waiting time $t-a(Y_t)$ of a delivery of grapes exceeds the acceptable limit, the winery can choose to either blend the grapes with a lower quality variety or discard that delivery altogether, thereby losing that potential income. Which option is selected depends on the amount of time that has elapsed since the grapes’ harvest. Respectively, the two possible options listed above translate into the following control-based changes: $Y_t$ changes the variety $v(Y_t)$ or, as the grapes are removed from the queue altogether, any action regarding them is removed from the set of available controls $\Gamma(t,X_{t-1})$.

Since wineries process a finite set of grape varieties and incoming trucks load weight can be modeled as a set containing a finite number of items, it is possible to create another discrete finite set of all possible truck types, describing incoming deliveries, and consisting of N possibilities $[v,l]\in\{[v_1,l_1],  \ldots, [v_N,l_N]\}$, where the vector $[v_i,l_i]$ specifies grape variety and load weight of the $i$th type, respectively. In order to emphasize truck dependence on the $i$th type, we employ the following shortened notation $[v(Y_t),l(Y_t)]=[v_i,l_i]=y_i$.

{\bf Set of new information $\mathbf{Z}$.}\label{sect_new_info}
At each time $t$, we can observe newly arrived trucks at the winery, which are  in turn represented by $Z_t$.
Each element $Z_t=\{Z_t^1,\ldots,Z_t^N\}$ is a set consisting of $N$ i.i.d. random variables denoted by a generic $Z_t^i$ ($i=1,\ldots,N$), with all $Z_t^i$ taking the form of \textit{Bernoulli random variables};
$$Z_t^i= B_t^i,\hspace{10pt}\textrm{where }B_t^i\sim \text{Bern}(p_t^i).$$
Each $B_t^i$ can assume one value in the set $\{0,1\}$, with probability
$$\mathbb{P}(B_t^i=1)=p_t^i,$$
where $B_t^i=1$ indicates that at least one truck of type $i$th arrives at time $t$.
Otherwise, when $B_t^i=0$, $Z_t^i=0$, no trucks of type $i$th arrive at time $t$.

{\bf Transition function $\mathbf{f}$.}\label{sect_trans_funct}
Recall that the following expresses the system state at time $t$:
$$X_t=f(t,X_{t-1}=x, Y_t=y, Z_t=\{z_1,\ldots,z_N\}),$$
where $f$ is a transition function and lower case symbols $x,y,z$ denote realizations of the corresponding  upper case random variables. For the winery possessing $K=1$ press it is specified as follows:
\begin{eqnarray*}
\lefteqn{\Big[v(X_t),l(X_t),s(X_t)\Big]=X_t=f\Big(t,x, y, \{z_1,\ldots,z_N\}\Big)}\\
&&=\sum_{i=1}^{N+1} \mathbbm{1}_{\big\{y=y_i\big\}}\mathbbm{1}_{\left\{\big(v(x)=v_i\hspace{3pt} {\rm or} \hspace{5pt}
v(x)=0\big)\hspace{5pt} {\rm and}\hspace{5pt} \big(C-l(x)\geq l_i\big)\right\}}\cdot\\
&&\Big[v_i\cdot\mathbbm{1}_{\left\{v_i\neq 0\right\}}+v(x)\cdot\mathbbm{1}_{\left\{v_i=0\right\}},\\
&& l(x)+l_i,\\
&& t\cdot\mathbbm{1}_{\left\{l(x)+l_i=C \hspace{3pt} {\rm and} \hspace{3pt} s(x)=0 \right\}}+s(x)\cdot\mathbbm{1}_{\{t-s(x)<TP \hspace{3pt} \rm and \hspace{3pt} s(x)\neq 0\}}
\Big].
\end{eqnarray*}
The function is constructed from a sum over $N+1$ possible options, referring to $N$ different truck types plus one symbolic truck $\mathbbm{O}$.  Truck $\mathbbm{O}$, defined as $y_0=[0,0]$, represents  the situation when no truck is assigned to the press being analyzed. Firstly, we have two indicators: $ \mathbbm{1}_{\{y=y_i\}}$, which extracts $[v_i,l_i]$ components from selected control $y$, and $\mathbbm{1}_{\left\{\big(v(x)=v_i\hspace{3pt} {\rm or} \hspace{5pt} v(x)=0\big)\hspace{5pt} {\rm and}\hspace{5pt} \big(C-l(x)\geq l_i\big)\right\}}$, which checks whether the grape variety the truck holds  matches the variety already inside the press being analyzed  or if the press is empty, under the assumption that this truck's load does not exceed the press's spare space. Only in such a context can truck $y$ be assigned to the press. Next, all components of $X_t$ are defined. The press changes variety to $v_i$ if $v_i$ does not equal 0. However, if the $\mathbbm{O}$ truck is assigned to the press, the variety must remain as it was preceding this assignment, since assigning the $\mathbbm{O}$ truck to a press does not change its state; instead, the weight is simply modified by adding selected truck load weight $l_i$ to the previous press state weight $l(x)$. If we fill the press $(l(x)+l_i=C)$, which has not yet started processing grapes $(s(x)=0)$, the processing start time $s(X_t)$ is assigned to the value $t$. If the press started processing $(s(x)\neq 0)$ and will not finish until $t$, $(t-s(x)<TP)$, we preserve the previous starting time $s(X_t)=s(x)$; otherwise this time is 0. Constant over-time factors $C$ and $TP$ remain the same.

{\bf Payoff $\mathbf{g}$.}\label{sect_payoff}
The optimized decision process can be termed optimal if it maximizes total profit gained from the wine production. Total profit consists of the sum of payoffs $g$, calculated at each time point $t$. In our model, every time the press begins a cycle, the manufacturer receives a payoff, defined as
\begin{equation*}
g(t,X_{t-1},Y_t)=
\begin{cases}
{\rm Price}(v(X_t))\cdot C,& s(X_t)=t, \\
0,& \text{otherwise}.
\end{cases}
\end{equation*}

The payoff is generated from state $X_{t-1}$, which is also affected by decision $Y_t$, and these return a new state $X_t$ according to the aforementioned transition function. $Price(\cdot)$ function defines the income received by the winery from pressing one unit of a certain grape variety. When the number of presses is strictly larger than $1$; that is, $K>1$, we sum the values of payoff function over all machines and optimize our decisions based on the sum of payoffs.

{\bf Value function $\mathbf{V}$.}\label{sect_value_funct}
All the elements of the Markov decision process described above allow us to define the objective, i.e., identifying the truck-assignment strategy $Y^*_{1\ldots T}$ that returns the highest expected total payoff from produced wine
$$ V_{Y^*_{1\ldots T}}(1,x_0)=V^*(1,x_0)=\sup_{Y_{1\ldots T}}V_{Y_{1\ldots T}}(1,x_0). $$

\section{Model calibration}\label{chapt_4}
In the previous sections, we introduced a theoretical model and provided a general interpretation of it, when applied to the grape-reception process. However, the proposed model is not straightforwardly applicable. It still requires calibration with respect to data and production-process knowledge obtained from representatives of a real-world winery.
In general, we know that the manufacturer, that graciously agreed to provide the information needed, owns six presses $\{X^1,\ldots,X^6\}$ of two different types. Recall from Section \ref{sect_model_def} that each press is described by a three-element vector comprised of variety, load weight, and processing start time and two additional invariant parameters -- total processing time and capacity. The type of press directly affects total processing time $TP$ and capacity $C$, and we distinguish the two press types as follows:
\begin{itemize}
\item Type I : \hspace{7pt} $C=$\SI{25}{\tonne},\hspace{2pt} $TP=$\SI{2}{\hour},
\item Type II : \hspace{3pt} $C=$\SI{50}{\tonne},\hspace{2pt} $TP=$\SI{4}{\hour}.
\end{itemize}
Out of six presses, four were of Type I, and the remaining two were of Type II.

The presses were used to process four different grape-varieties, which could not be blended. Although each grape-variety  generated wine having a different price, the winery’s representatives did not disclose information pricing information. We assigned the grape-varieties the numbers 1 through 4, which we straightforwardly assumed to be the payoff obtained from processing 1 metric ton of appropriate variety of grapes; i.e., the payoff from processing 1 metric ton of grape-variety 1 is assumed to be 1, the payoff from processing 1 metric ton of grape-variety 2 is 2, etc.

As already mentioned in Sections \ref{chapt_1} and \ref{sect_model_def}, once harvested, grapes deteriorate with time, and the winery distinguishes two phases of deterioration. The first phase begins after 2\si{\hour} of waiting in the queue to be unloaded. Grape quality then drops to the level of those grapes used to manufacture the cheapest variety of wine and can therefore be blended only with the types of grapes normally employed to make this wine type. The second phase, which begins after 4\si{\hour} of waiting in the queue, sees quality drop to the level which makes the grapes unacceptable for wine production altogether and the delivery has to be disposed. In this case truck after
waiting 4 hours in a queue is removed from the queue.

Further calibration, such as that needed to approximating probability distribution of delivery times, required historical queue data. The data provided by the winery, which was not  consistent in type of information provided, consisted of two datasets. The first dataset contained the total weight of grapes that arrived at the manufacturer on a given date, grouped by grape-variety. The second dataset, however, contained information about the load weight of the grapes each truck contained and its arrival time but, unfortunately, did not include grape variety. On the basis of the second dataset, we were able to discretize the time. The dataset was aggregated to the form presented in Table \ref{tab_opening}, identifying hours of each day’s first and last delivery as well as the average time between deliveries. Taking into account that deliveries began arriving from 8:30  and continued throughout the day to about 00:30, we divided one day into $34$ half-hour intervals from 8:30 to 01:30, with an additional one hour added for managing deliveries from the last lot. Such a split forces the system to deal with approximately $4$ trucks per half-hour interval. Given  that the maximum number of available machines is $6$, the proposed split is reasonable.

\begin{table}[h]
\centering
\begin{tabular}{|c|c|c|c|}
 \hline
 \rowcolor{lightgray}
 \textbf{Day of week} & \textbf{First delivery} & \textbf{Last delivery} & \textbf{Average time between deliveries [\si{min}]}\\
  \hline
  Mon. & 08:46 & 00:42  &  7.41\\
  Tue. & 08:38 & 00:17  &  6.84\\
  Wed. & 08:37 & 23:05  &  7.23\\
  Thu. & 08:45 & 23:00  &  7.77\\
  Fri. & 08:27 & 23:55  &  10.09\\
  Sat. & 08:42 & 23:56  &  5.78\\
  \hline
\end{tabular}
\vspace{0.5cm}
\caption{Average time between deliveries for each working day of a week}
\label{tab_opening}
\end{table}

On the basis of the second dataset, we created frequency plots of delivery loads weight for each day, and these are presented in Figure \ref{fig:weight_hist}. As can be seen, load weight ranged from 0 up to \SI{25}{\tonne} per delivery. We discretized the loads weight to the \SI{5}{\tonne}  intervals, obtaining a set $\{5,10,15,20,25\}$ of all possibilities. The discretization process assigned \SI{5}{\tonne} to a load whose weight fell within the $(0,5]$ interval, \SI{10}{\tonne} to a load with a weight falling within the $(5,10]$ interval, etc.
Recall from Section \ref{sect_model_def} that the combination of a discrete set of loads weight and a discrete set of grape varieties creates a discrete set of all possible truck type  designations $\{[v_i,l_i]: i=1,\ldots,N\}$, with $N=20$ possibilities in this case, resulting from combining four possible grape variety's types and five load weight possibilities.

\begin{knitrout}
\definecolor{shadecolor}{rgb}{0.969, 0.969, 0.969}\color{fgcolor}
\begin{figure}[h]
{\centering \includegraphics[width=\maxwidth]{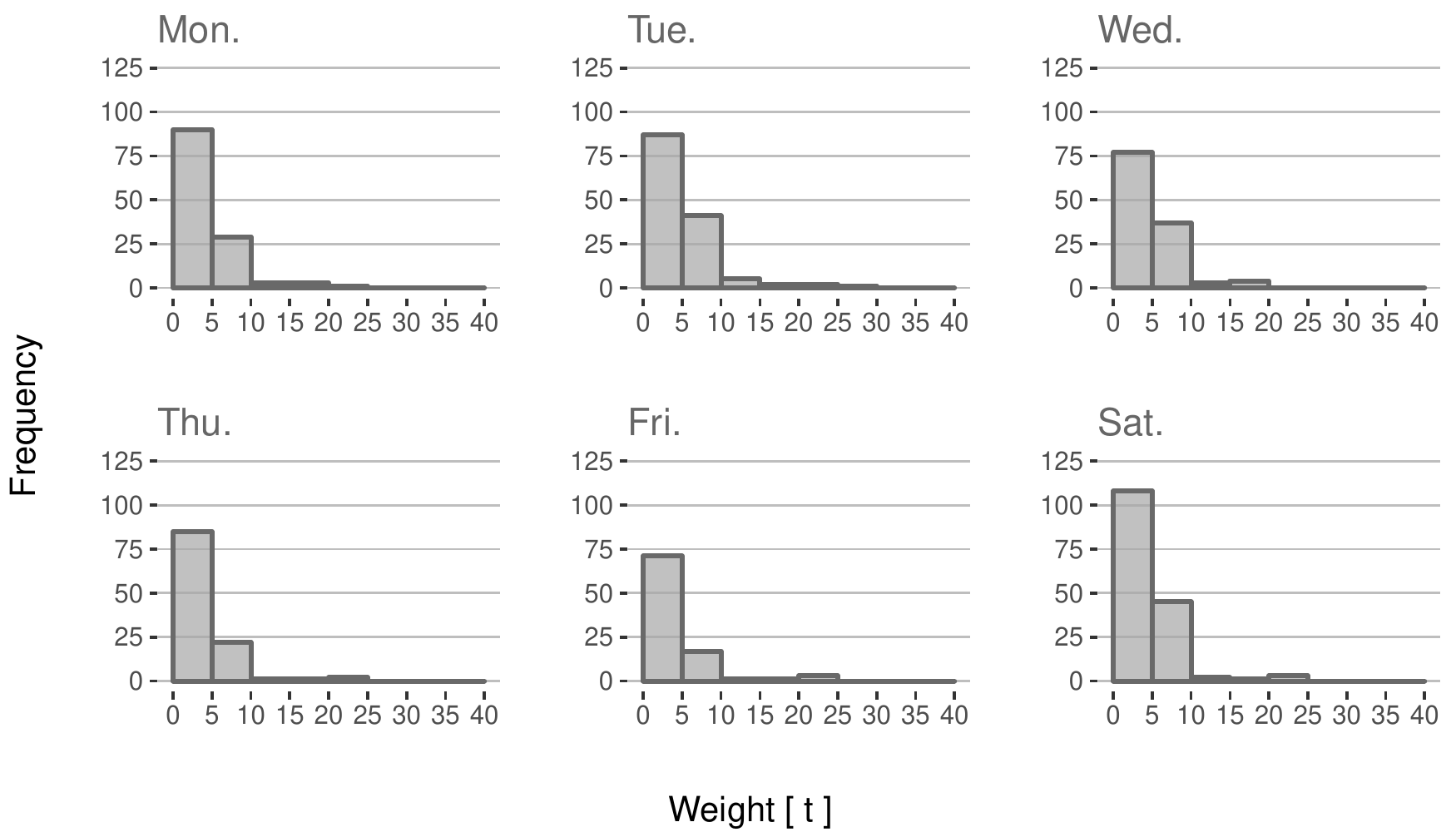}}
\caption{Frequency plots of delivery loads weight in metric tonnes for each working day of a week}\label{fig:weight_hist}
\end{figure}
\end{knitrout}
To summarize  all information presented so far, we define the following domains:
\begin{itemize}
\item \textbf{time} \hspace{21pt} $t\in \{0,1,\ldots,33\}$,
\item \textbf{grape variety} \hspace{4pt} $v_i\in \{1,2,3,4\}$, where $i=1,\ldots,20$,
\item \textbf{weight} \hspace{21pt} $l_i\in\{5,10,15,20,25\}$ in metric tonnes, where $i=1,\ldots,20$,
\item \textbf{capacity} $C\in\{25,50\}$ in metric tonnes,
\item \textbf{total processing time} $TP\in \{4,8\}$ in time intervals $t$,
\item \textbf{load's variety-weight combination} $i\in \{1,2,\ldots,20\}$.
\end{itemize}
The other time-related variables \textbf{arrival time} and \textbf{processing start time} are defined in the same way as \textbf{time}. The domains listed above allow us to completely define two model's components: the state of the stochastic dynamical system ($X$) and the set of model's controls ($Y$).

\begin{figure}[H]
\includegraphics[width=0.7\textwidth,height=0.25\textheight]{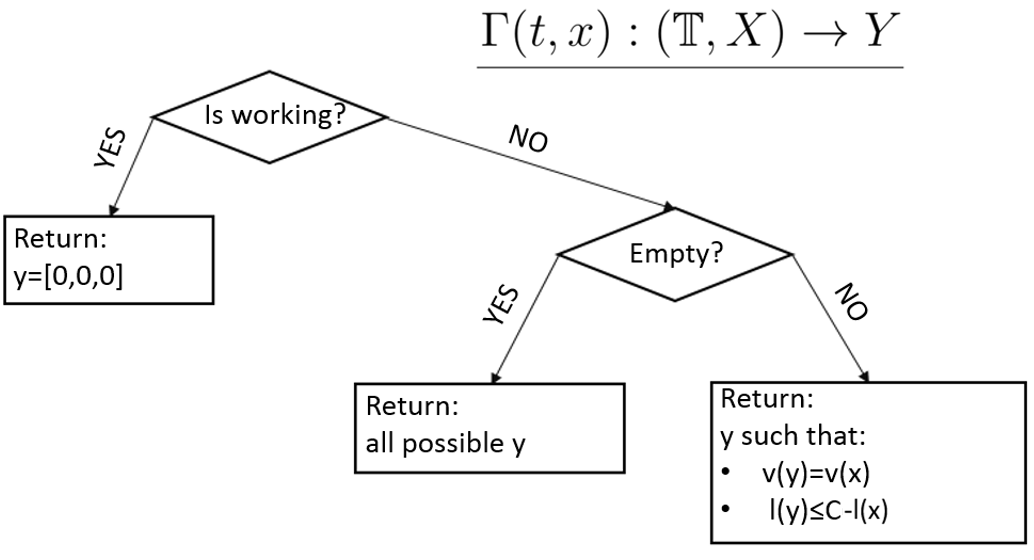}
\centering
\caption{Decision tree for $\Gamma(t,X)$}
\label{fig_dec_tree}
\end{figure}

The next component we address is $\Gamma$, the set of available decisions given the system’s current state, can most comprehensibly be designated by the decision tree presented in Figure \ref{fig_dec_tree}. Were a press to be blocked, filling it with any delivery would be impossible, in which case the decision tree would return the symbolic truck type $\mathbbm{O}$. If a press is not operating, the decision tree poses the next question, which concerns filling the press. For an empty press, the decision tree returns all possible truck types, where truck type is represented by one of 20 variety-weight combinations; otherwise, it returns only those trucks having the same grape variety as that of the grapes the press currently contains provided their loads do not exceed the press's spare space. As presses are of two types, $\Gamma$ is also defined twice, i.e., Type I and Type II, respectively.

Recall from Section \ref{sect_model_def} that the set of new information $Z_t$ consists of $N=20$ Bernoulli random variables, each indicating the presence of the truck types among the deliveries occurring during time interval $t$. That is,
\begin{equation}\label{pi_indep}
\mathbb{P}\left(Z_t^i=1\right)=\mathbb{P}(y_i)=\mathbb{P}\left([v_i,l_i]\right)=\mathbb{P}(v_i)\cdot \mathbb{P}(l_i)=p_i, \hspace{10pt}\textrm{where } i=1,\ldots,N.
\end{equation}
At the moment, we assume the following. During one time interval, the manufacturer obtains the load of at most one truck; the variety of the grapes being delivered is independent of the load weight  they comprise; and the probability of receiving each truck type (i.e., denoted by a combination of truck’s load weight and variety of grapes contained) is constant over time. In such a context, approximating the probability of a truck’s load being of one grape variety based on that variety’s frequency within the first dataset, the probability of a discretized load weight presence  within the second dataset, and use of the independence property to calculate the probability of a certain truck type  by marginal probabilities multiplication is possible. Figure \ref{histogram_tu_22} shows the example data from one day.

\begin{knitrout}
\definecolor{shadecolor}{rgb}{0.969, 0.969, 0.969}\color{fgcolor}
\begin{figure}[h]
{\centering \includegraphics[width=\maxwidth]{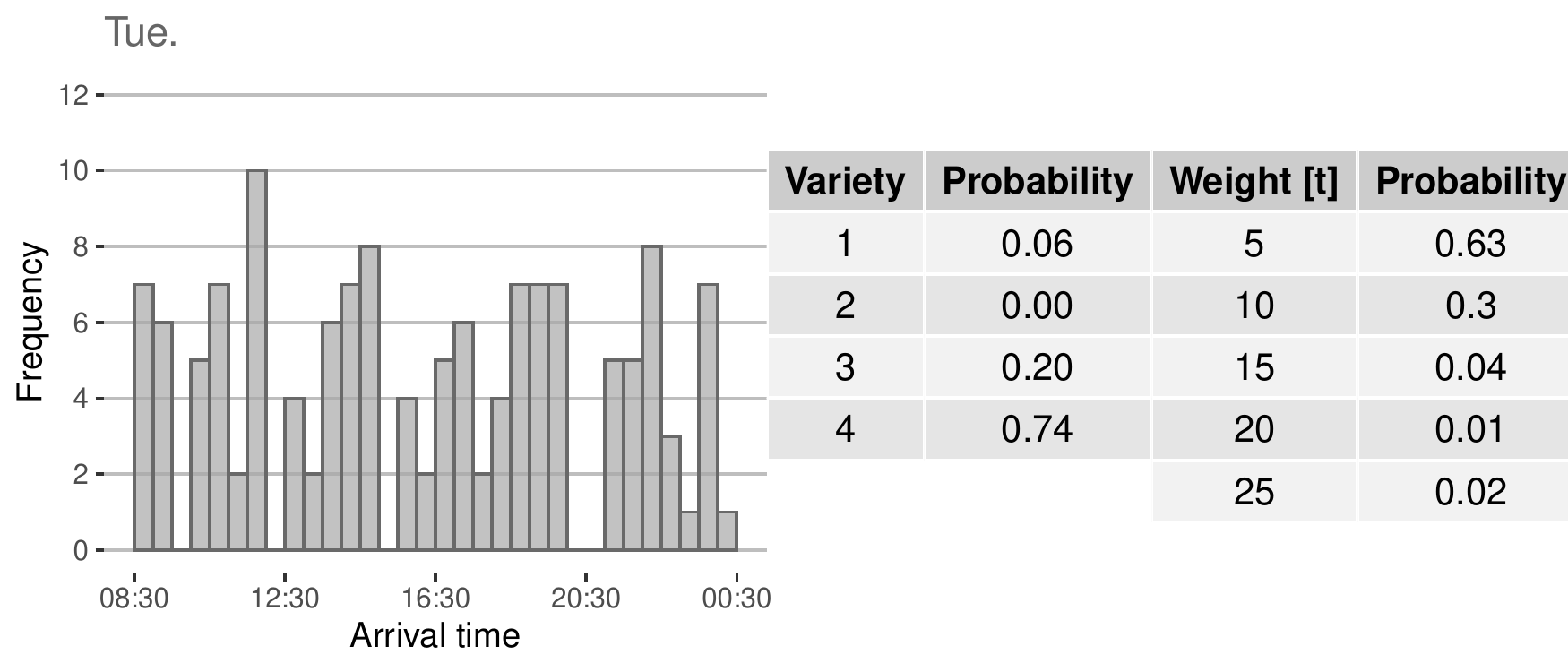}}
\caption{Tuesday frequency plot of deliveries with marginal probabilities for variety and weight \label{histogram_tu_22}}
\end{figure}
\end{knitrout}

\begin{knitrout}
\definecolor{shadecolor}{rgb}{0.969, 0.969, 0.969}\color{fgcolor}
\begin{figure}[h]
{\centering \includegraphics[width=\maxwidth]{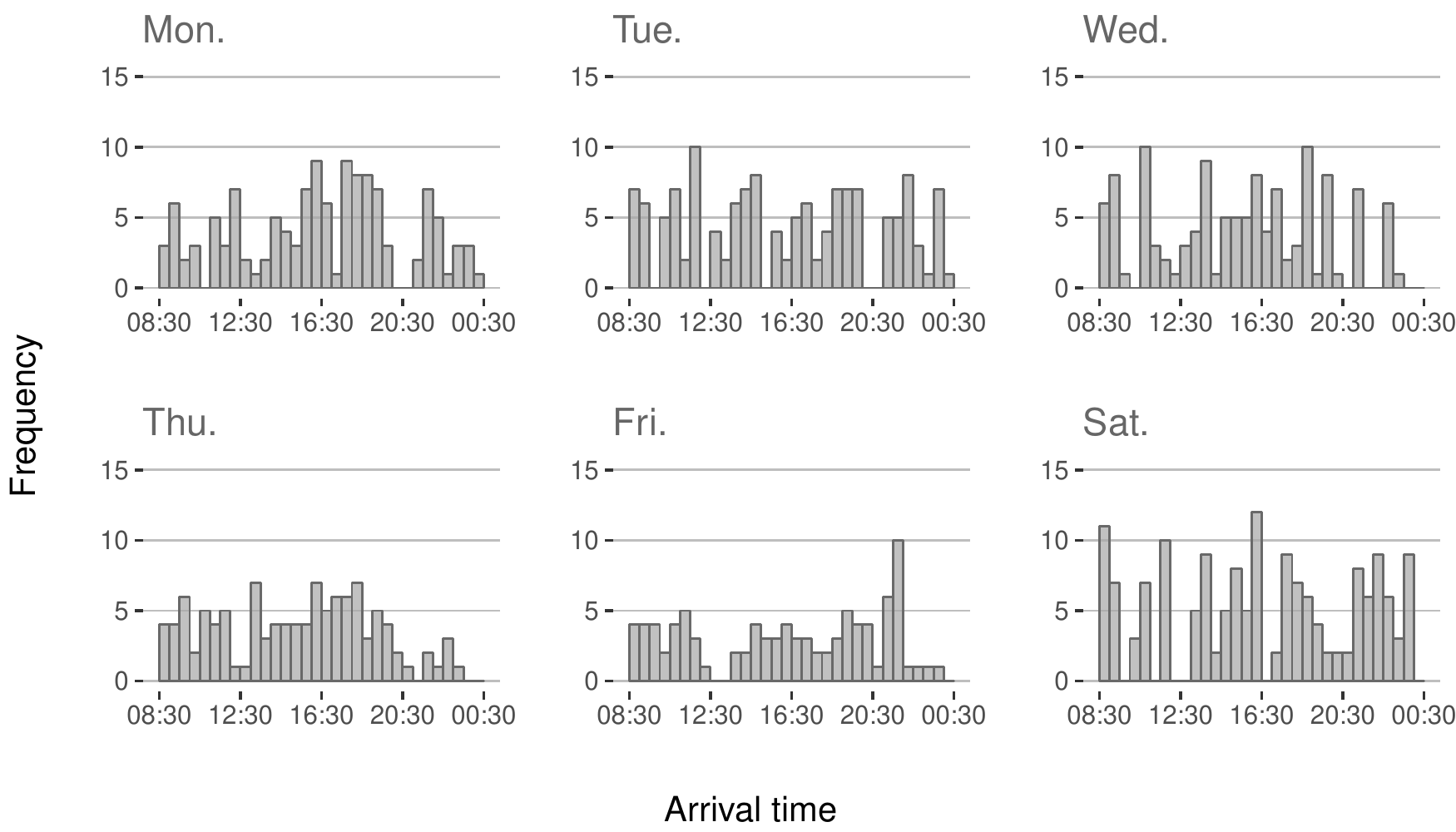}}
\caption{Frequency plots of deliveries amount for each working day of a week}\label{fig:deliveries_hist}
\end{figure}
\end{knitrout}

As can be seen in the frequency plots shown in Figure \ref{histogram_tu_22} and Figure \ref{fig:deliveries_hist}, the assumption of at most one truck arrival during each half-hour interval does not reflect reality. $Z_t$ distribution must be complemented with the information about number of deliveries in each interval. At first glance, although the distribution appears to be somewhat random, a specific ‘gap-structure’ is noticeable for every day (around 10 AM, 12 AM, and 8 PM). To check whether this pattern is repetitive, we generated stripe plots (see Figure \ref{fig:poiss_hist}), where each stripe marks a delivery arrival on a continuous time scale and added to that a line showing the average number of deliveries in a 30-minute time interval, marked with red data points and placed in the middle of the represented interval.

\begin{knitrout}
\definecolor{shadecolor}{rgb}{0.969, 0.969, 0.969}\color{fgcolor}
\begin{figure}[h]
{\centering \includegraphics[width=0.85\maxwidth]{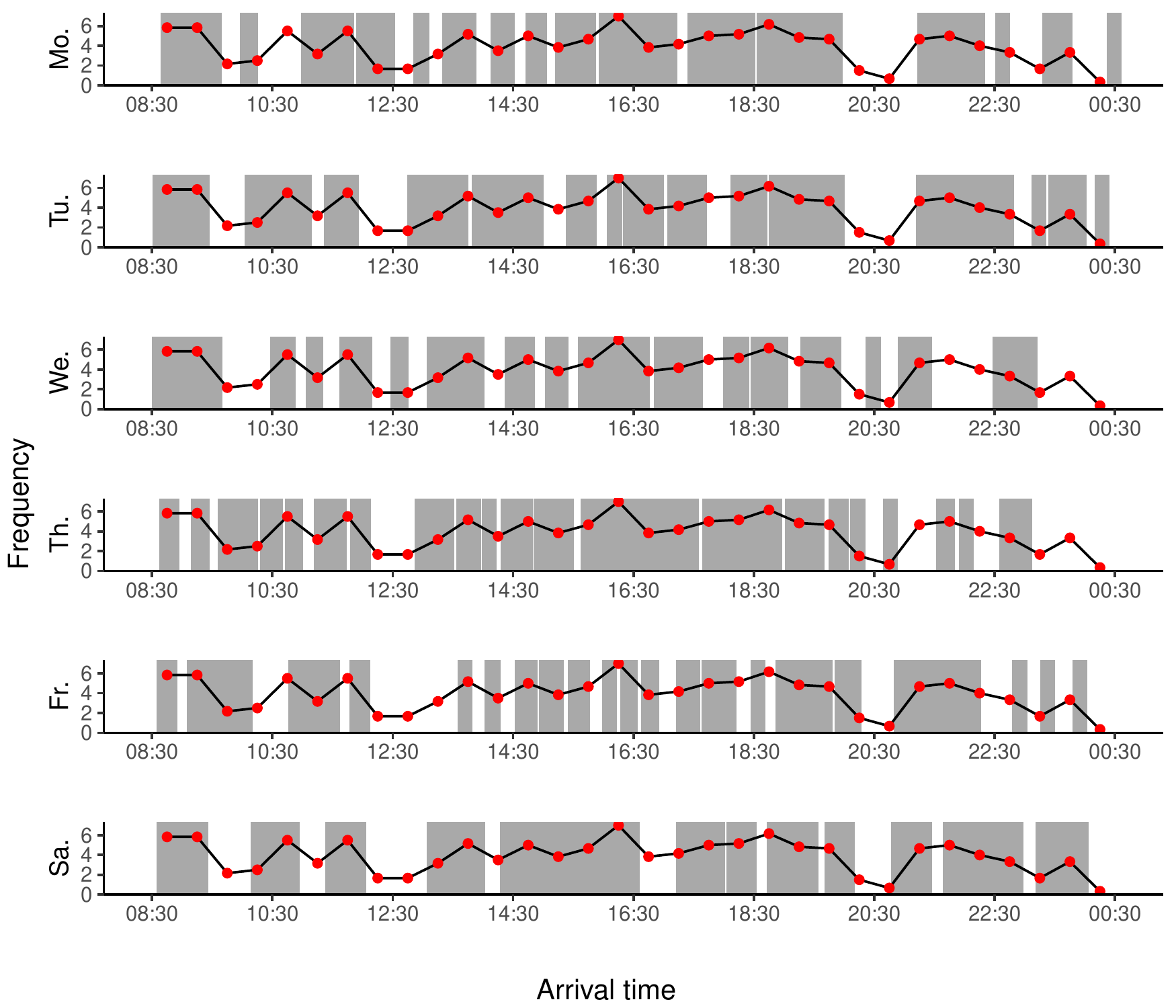}}
\caption{Stripe plots showing delivery amounts for each working day of a week}\label{fig:poiss_hist}
\end{figure}
\end{knitrout}

Figure \ref{fig:poiss_hist} shows a specific, repetitive delivery pattern over time, suggesting that the data may not be completely random. Since the most common approach to modeling arrivals is the doubly stochastic Poisson process, we applied that approach to this process, with a different $\lambda$ parameter for each time interval $t$. We approximated this parameter with the mean number of deliveries within each time interval, across the days, and employ red dots to represent this quantity in Figure \ref{fig:poiss_hist}. We further represent the delivery counting process in interval $t$ by variable $D_t$, which in the first set of assumptions, i.e., one truck per one time interval, equaled at most 1. When $D_t\sim {\rm Poiss}(\lambda_t)$, then

\begin{equation}\label{poiss_D}
\mathbb{P}(D_t=k)=\frac{\lambda_t^k e^{-\lambda_t}}{k!}, \hspace{10pt}\textrm{where } k=0,1,\ldots
\end{equation}

Combining (\ref{pi_indep}), (\ref{poiss_D}), and the $e^x$ series expansion yields the probability that the $i$-th truck type is present among the deliveries occurring in interval $t$
\begin{align}\label{PDF_Z}
\mathbb{P}\left(Z_t^i=1\right)&=\sum_{k=0}^\infty\mathbb{P}(D_t=k)\cdot\mathbb{P}\left(Z_t^i=1|D_t=k\right)=\sum_{k=0}^\infty \frac{\lambda_t^k e^{-\lambda_t}}{k!} \cdot \left(1-(1-p_i)^k\right)\\
&= 1-e^{-\lambda_t p_i}.\nonumber
\end{align}
When interpreting the probability given above, we observe that the larger $\lambda_t$ or $p_i$ is, the greater is the likelihood that $Z_t^i$ equals 1.

\begin{knitrout}
\definecolor{shadecolor}{rgb}{0.969, 0.969, 0.969}\color{fgcolor}
\begin{figure}[H]
{\centering \includegraphics[width=0.85\maxwidth]{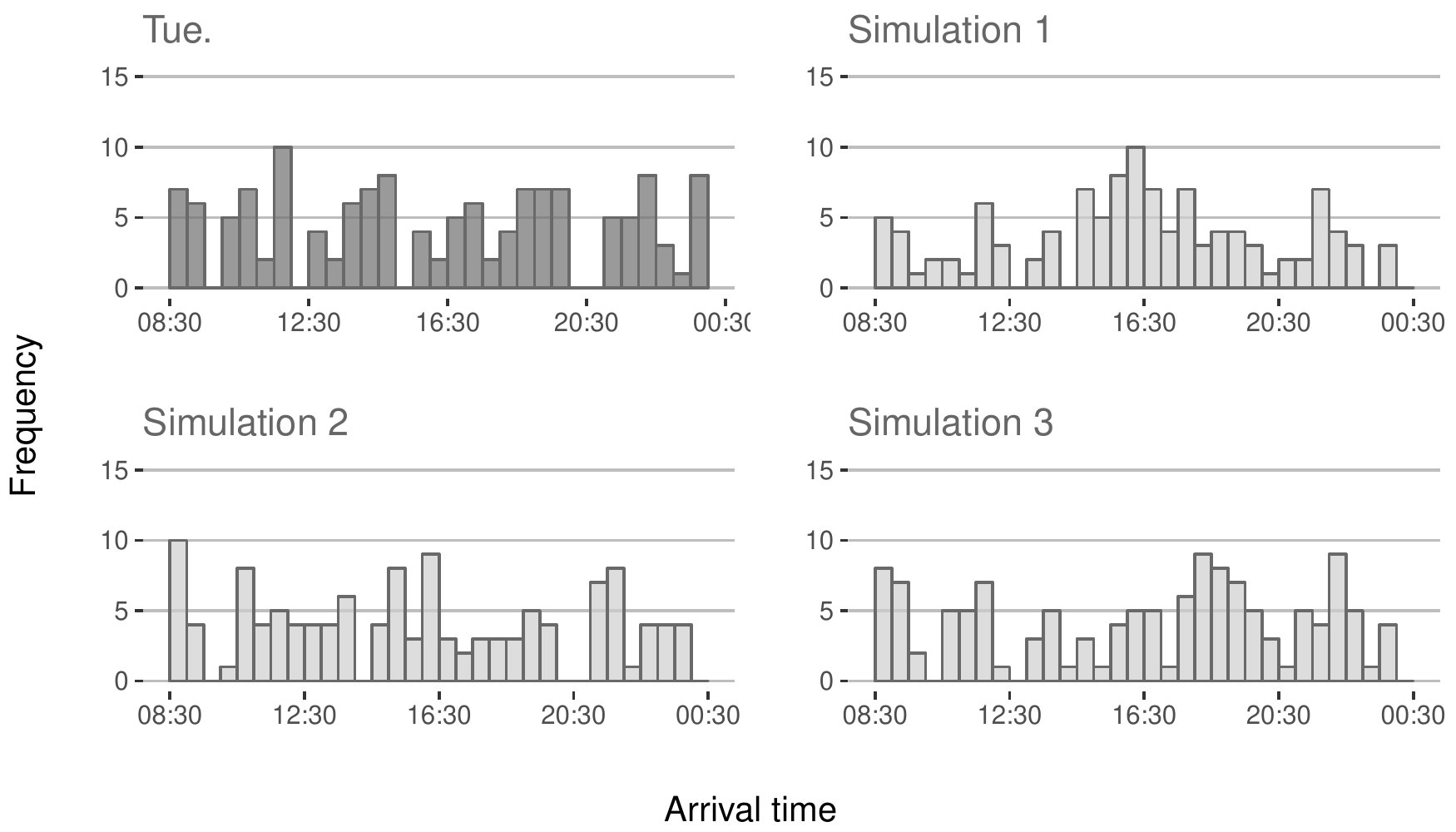}}
\caption{Daily  frequency plots of simulated deliveries arranged with simulation-base Tuesday data\label{arrivals_simulated}}
\end{figure}
\end{knitrout}

Figure \ref{arrivals_simulated} presents frequency plots of the original Tuesday queue and simulated queues, based on the arrivals being a doubly stochastic Poisson process.

Recall that the winery’s presses can be of two types, implying that the value function $V$, which is reliant on a press's state, should be calculated twice, once for each version. This difference in state is generated by a $\Gamma$ function that filters controls $Y$ for the value function; that is, a Type II press can accommodate a larger number of controls because its capacity (\SI{50}{\tonne}) is larger than that of a Type I press (\SI{25}{\tonne}).

\section{Algorithm implementation}\label{chapt_5}
With the theory  underlying the model known and the model’s components having been described and defined, the last step before obtaining a maximizing strategy was implementation of the optimization algorithm. At this point, two crucial questions need to be answered:  What is the objective of the optimization problem? Over which arguments is the optimization performed?
Due to the fact that the manufacture possesses several presses but operates a single queue common to all the presses, optimization has to be performed for the sum of value function $V$ outputs for presses $X^1,\ldots,X^6$, over possible delivery permutations. Such an approach adheres to the constraint that each delivery can be assigned only once.

The other important constraint is grape deterioration over time, which the algorithm incorporates by including an overall deterioration cost $c$ in value function $V$. Mathematically, this can be obtained by subtracting from value function the deterioration expenses resulting from each decision. Among the expenses are variety's degradation costs
$${\rm c}_{\rm deg}(Q_{\rm deg})=\sum_{i=1}^{M_{\rm deg}} \left(v(Q_{{\rm deg}_i})-1\right)\cdot l(Q_{{\rm deg}_i})\hspace{10pt}$$
and rejection costs
$${\rm c}_{\rm rej}(Q_{\rm rej})=\sum_{i=1}^{M_{\rm rej}} l(Q_{{\rm rej}_i}).$$
Here $Q_{\rm deg}$ and $Q_{\rm rej}$ define sets of trucks whose content is classified for degradation or rejection, respectively, where each truck is represented so as the control $Y$. The $M_{\rm deg}$ and $M_{\rm rej}$ define the numbers of such trucks.

The manufacturer allows splitting of deliveries, i.e.,\ accepting only \SI{5}{\tonne} from a \SI{20}{\tonne} delivery, and combining deliveries, i.e.,\ assigning \SI{5}{\tonne} from one truck and \SI{5}{\tonne} from another truck to a single press. The only limitations on such decisions are to keep them consistent with the $\Gamma$ function and not to exceed the threshold of distributing at most \SI{75}{\tonne} of grapes during one time interval. The latter limitation arises from the winery’s workflow capacity, including grape unloading and preprocessing prior to pressing.

Applying the rules given above to the pre-specified model results in an algorithm that quite closely reflects reality.
As an initial input, the algorithm requires parameters for the distribution of expected deliveries, and these parameters are passed to the algorithm in three lists. The first list defines the variety's ratio among the total amount of grapes being delivered, i.e., the probability of each grape variety’s arrival. The second list defines the grape-weight-per-truck ratio in the total amount of grape deliveries. Recall from Section \ref{chapt_4} that the load weights are discretized based on five intervals; the second list returns the probability for each of these. The third list contains the expected arrival intensity of trucks per each time interval.

This initial input is then transformed into truck type distribution $Z$, according to the probability function given in Equation (\ref{PDF_Z}). The output gives the probabilities that a time interval contains the truck types,  which, in turn, are used to calculate the total expected payoff resulting from a given state, by the means of Bellman equation (\ref{recursion_eq}). This total payoff is the summarized payoff obtained at the end of each day. The algorithm returns a three-dimensional list of tables: one table per each time interval, and every table contains the maximal expected total payoffs at the end of the day for every possible state at the corresponding time interval. The list is generated in two versions, one for each type of press, and these  outputs serve as the basis for further choice within the profit-maximizing strategy.

Now that we have described all required inputs to the model, we can proceed to describe the model’s decision-making procedure, which is divided into two functions. The first function returns information related to maximizing press fill, i.e., how many metric tonnes of which grape variety should be assigned to each press. As an input, it takes a queue containing information about each truck load variety, weight and arrival time; the presses' states; the current time  interval; and the maximal expected total payoff for all states possible at this time interval (denoted by \textit{V\_T}). At first, the function checks to see if the queue is empty, and, if it is, the function generates the decision not to fill the presses. Otherwise, if the queue does contain one or more trucks, the function checks to see how the presses could be filled, according to their present states ($GP$), and defines a variable describing the total weight of each grape variety  ($Queue\_aggr$) that could be unloaded. In the next step, the possibilities are merged with the payoff they would generate in order to enable assessment of each possibility. These possibilities are then compared with the truck types and numbers comprising the queue, so as to choose a valid fill decision. In addition to income, the algorithm observes loss generated by degraded or wasted truck loads remaining following each decision. At the end, the function returns all decisions that maximize profit, which is an income reduced by a loss.

\begin{algorithm}[H]
\caption{Chooses maximizing presses' fill}\label{fill_decisions}
\begin{algorithmic}[1]
\Function{FILL\_DECISIONS}{$Presses,Queue,time, \textit{V\_T}$}
\If{$\text{Queue is empty}$}
\State $\textit{max\_Y} \gets \mathbbm{O}$ \Comment{do not fill the presses}
\Else
\State $\textit{GP} \gets \textit{Gamma(Presses)}$ \Comment{all possibilities to assign for each press}
\State $\textit{Queue\_aggr} \gets \text{Queue aggregated by variety}$
\State $\textit{GP\_payoff} \gets merge(GP,\textit{V\_T})$
\State $\textit{Y\_income} \gets \textit{adjust(GP\_payoff,Queue\_aggr)}$
\State $\textit{Y\_loss} \gets \textit{loss(Y\_income,Queue,time)}$
\State $\textit{Y\_profit} \gets \textit{Y\_income-Y\_loss}$
\State $\textit{max\_Y} \gets \textit{max(Y\_profit)}$
\EndIf
\State $\textbf{return } \textit{max\_Y}$
\EndFunction
\end{algorithmic}
\end{algorithm}

The second function returns a strategy consisting of instructions on how the grapes should be unloaded from the trucks. The general rule states that the trucks should be unloaded in chronological order, from the earliest to the latest arrived truck. The input is identical as for the \textit{FILL\_DECISIONS} function, since it is executed from the inside of the second function. At the beginning, the function initiates a variable for storing the strategy. Then, by the use of a \textit{FILL\_DECISIONS} function, it generates the maximizing possibilities of filling the presses. If among the possibilities at least one differs from `do not fill the presses', the function randomly picks one of the maximizing strategies and fills the presses using the oldest possible trucks. At the end, it returns a maximizing strategy in the form of a table containing information about the press filled, the id of the truck used to fill it, the truck arrival time, the truck variety, and the amount of the truck load utilized. In case of only one maximizing possibility `do not fill the presses', the function returns an empty table.

\begin{algorithm}[H]
\caption{Assigns trucks to presses}\label{run_model}
\begin{algorithmic}[1]
\Function{RUN\_MODEL}{$Presses,Queue,time,\textit{V\_T}$}
\State $\textit{max\_strategy} \gets \text{initiation of variable storing trucks usage}$
\State $\textit{max\_fill} \gets \textit{FILL\_DECISIONS(Press,Queue,time,V\_T)}$
\If{$\textit{any(max\_fill)} \neq \mathbbm{O}$}
\State $fill\_dec \gets \text{randomly chosen one of the maximizing strategies}$ 
\For{$P\text{ in } Presses$}
\State $\textit{[chrono\_trucks,Queue]} \gets \textit{Chronological\_truck\_use(Queue,max\_fill[P])} $
\State $\textit{max\_strategy} \gets \textit{append(max\_strategy,chrono\_trucks)}$
\EndFor
\EndIf
\State $\textbf{return } \textit{max\_strategy}$
\EndFunction
\end{algorithmic}
\end{algorithm}

It is worthwhile to notice that, for superfluous time intervals $(t\in\{T-1,T\})$, further maximal expected total payoff is set to 0, a decision deriving straightforwardly from the fact that $\mathbb{P}(Z_{t}^i=1)=0$, for each $i=1,\ldots,20$ and $t\geq T-1$, since deliveries arrive to the manufacturer from $t=0$ to $t=T-2=31$. In such cases, the optimal value would be based only on the current fill evaluation.

The proposed algorithm was implemented in the R programming language (version 3.6.1), and all computations were run on a computer having the Intel Core i5-3340M processor running at 2.7 GHz with 12GB of RAM under 64-bit Windows 10. Approximately 1 minute and 10 seconds were required to generate a table with maximal expected total payoffs for all states in two versions, and approximately 2 minutes were needed to select a maximizing strategy at a single time interval. All scripts used in this article are accessible from \cite{link}.

\section{Results}\label{chapt_6}
\begin{figure}[h]
\includegraphics[width=\maxwidth]{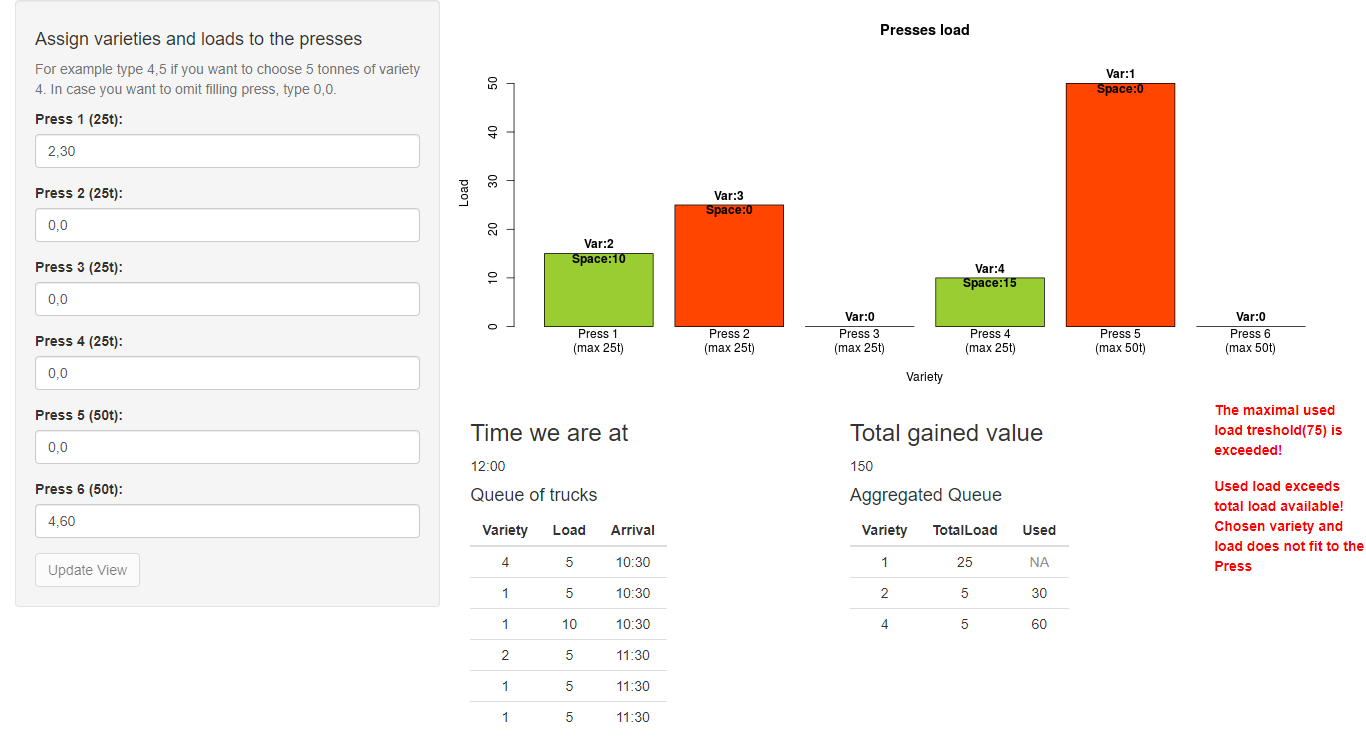}
\centering
\caption{Interface snapshot with error handling}\label{figinterface}
\end{figure}

In order to assess the algorithm, we compared the results the algorithm generated with those based on the method currently used by the manufacturer. The winery conducts the delivery-assignment process manually, with no software support. The manager assigns deliveries to the presses according to his expert knowledge and an approximated delivery schedule. To reflect this process, we created a graphical user interface, in which a user undertook decisions so as the manager in the winery did. The application was implemented in the R programming language using the 'Shiny 1.4.0' package and deployed  on {\it shinyapps.io}. The results generated by the users were then passed to an external file on hosting service and further analyzed. Figure \ref{figinterface} presents a snapshot of the application. The program can be found under the \cite{link2}.

\begin{figure}[!htb]
\includegraphics[width=0.95\maxwidth]{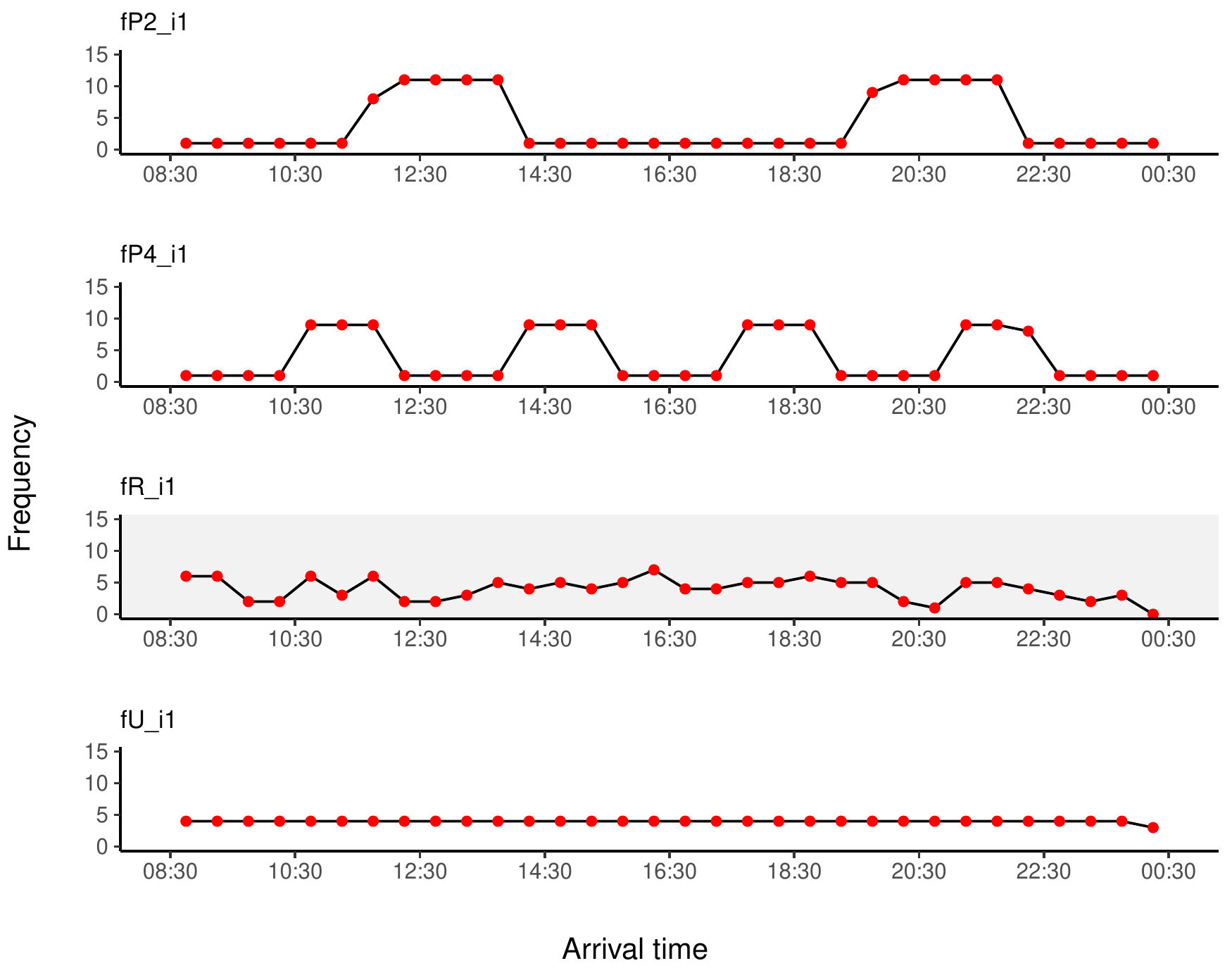}
\centering
\caption{Deliveries frequency \textbf{f} parameters, preserving real intensity}
\label{fig:freq_lambda}
\end{figure}

\begin{figure}[!htb]
\includegraphics[width=0.95\maxwidth]{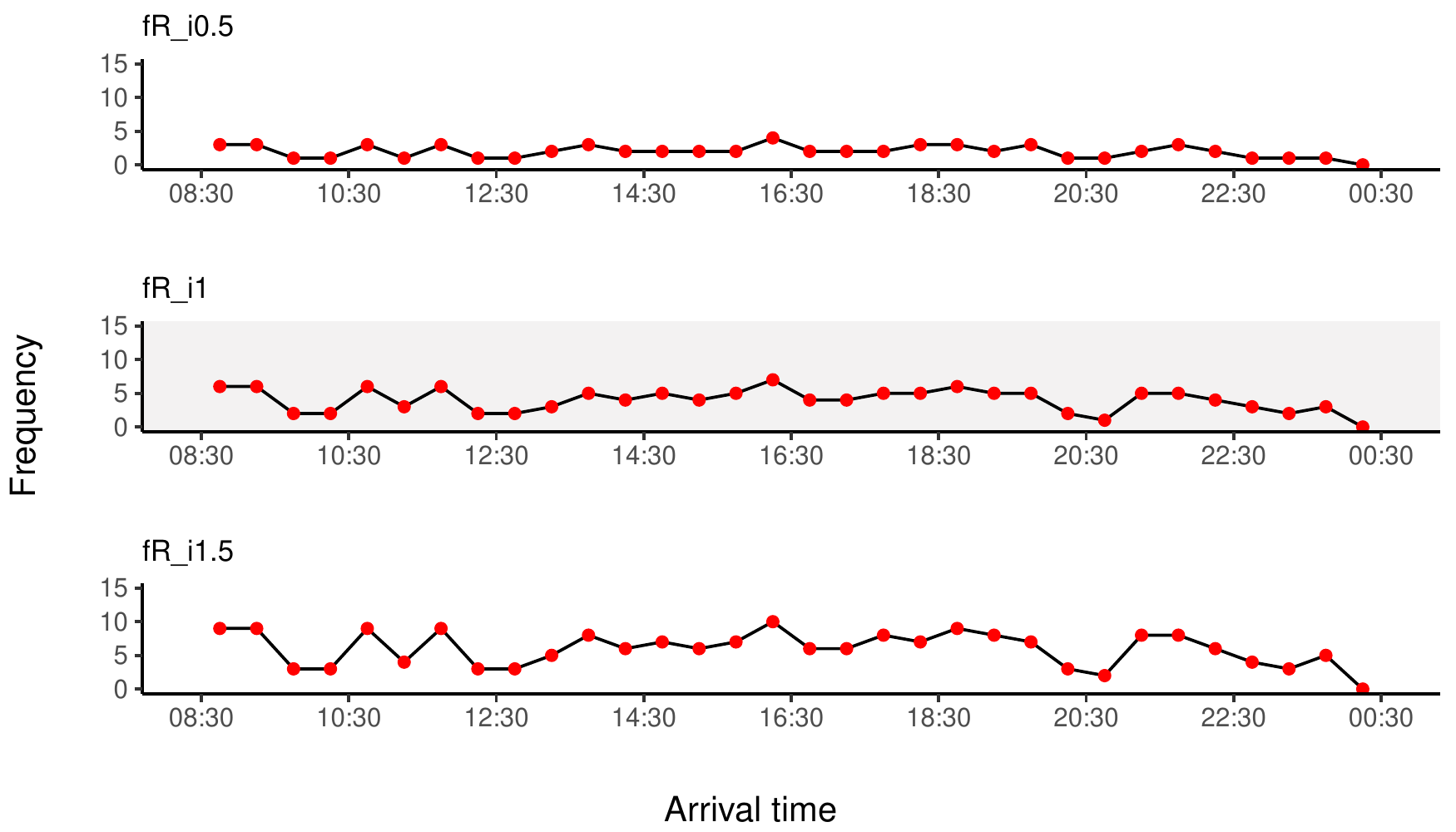}
\centering
\caption{Delivery intensity \textbf{i} parameters, preserving the real `shape' of frequency }\label{fig:int_lambda}
\end{figure}

While testing the algorithm, we identified three aspects of the decision process  that we wanted to research further in order to observe their influence on the results. The aspects are as follows:
\begin{itemize}
\item \textbf{Variety dominance v} defines which variety or grape varieties are over-represented in the testing scenario. Additionally, we inspected real variety's distributions from the Tuesday data, presented in Figure \ref{histogram_tu_22}, which is marked as vR. The specific parameters tested are shown in Table \ref{tab:variety_dominance_params}.\item \textbf{Deliveries frequency f} describes the distribution of deliveries over time. Among the testing scenarios, we observed one instance of two high but rare peaks (fP2), another with four lower but more frequently occurring peaks (fP4), and a uniform instance (fU) with a steadily occurring, equal number of deliveries. Additionally we included real-world frequency based on Tuesday data  (fR). Important to notice is that the total number of deliveries defined by this parameter is the same for each instance. The parameters are presented in Figure \ref{fig:freq_lambda}.
\item \textbf{Deliveries intensity i} describes the scaling of the total number of deliveries, preserving the `shape' of the frequency, based on real-world Tuesday data. The investigated test cases included a 0.5 times scaling of the total number of deliveries (i0.5), no scaling (i1) of this total number, and a 1.5 times scaling of the total number of deliveries (i1.5). Figure \ref{fig:int_lambda} presents these parameters.
\end{itemize}
All the simulations were described with the name v?\_f?\_i?, which explicitly defined the scenario being tested. During tests, we manipulated one parameter, preserving Tuesday’s real-world values for the two remaining parameters. The results were then averaged over four different realizations of the queue, satisfying the v,f,i parameters, which, for 21 different testing scenarios, yielded 84 simulations. The averaged results of these simulations can be found in Appendices Section \ref{apx:DA_sim_results}.

In general, the algorithm generated  \textbf{better} results than manual assignment. Specifically, the results generated by the algorithm were on average 5.69\% higher than those obtained through manual assignment, meaning that, for the average daily profits equal to 2364 generated through manual assignment, applying the proposed algorithm would instead increase total daily profit by 135. Recall from Section \ref{chapt_4} that the payoff is measured as the number of pressed metric tonnes of grapes multiplied by the variety's identifier. If we assume that the income from producing one bottle of wine approximately equals its variety's identifier multiplied by 1{\small \euro}, and approximately 1\si{\kilogram} of grapes is needed to produce one bottle of wine (see \cite{vinho_verde,price_wine}), then applying the algorithm would increase daily income by 135~000{\small \euro}.

\begin{table}[!htb]
\centering
\begin{tabular}{|l|cccc|}
  \hline
\rowcolor{lightgray} \textbf{Description} & \textbf{1} & \textbf{2} & \textbf{3} & \textbf{4} \\
  \hline
v1 & 0.7 & 0.1 & 0.1 & 0.1 \\
  v12 & 0.4 & 0.4 & 0.1 & 0.1 \\
  v123 & 0.3 & 0.3 & 0.3 & 0.1 \\
  v124 & 0.3 & 0.3 & 0.1 & 0.3 \\
  v13 & 0.4 & 0.1 & 0.4 & 0.1 \\
  v134 & 0.3 & 0.1 & 0.3 & 0.3 \\
  v14 & 0.4 & 0.1 & 0.1 & 0.4 \\
  v2 & 0.1 & 0.7 & 0.1 & 0.1 \\
  v23 & 0.1 & 0.4 & 0.4 & 0.1 \\
  v234 & 0.1 & 0.3 & 0.3 & 0.3 \\
  v24 & 0.1 & 0.4 & 0.1 & 0.4 \\
  v3 & 0.1 & 0.1 & 0.7 & 0.1 \\
  v34 & 0.1 & 0.1 & 0.4 & 0.4 \\
  v4 & 0.1 & 0.1 & 0.1 & 0.7 \\
  vU & 0.25 & 0.25 & 0.25 & 0.25 \\
  vR & 0.06 & 0 & 0.2 & 0.74\\
   \hline
\end{tabular}
\vspace{0.5cm}
\caption{Variety dominance  \textbf{v} parameters. Digits indicate the dominating varieties}
\label{tab:variety_dominance_params}
\end{table}

While preserving real-world frequency and intensity, the algorithm performed on average 6.31\% better than the manual-assignment method. The greatest predominance of the proposed algorithm over the manual method, 11.73\%, was obtained for a uniform distribution of grape varieties (vU), whereas the worst results, 0.44\% lower than the manual method, were generated when variety 2 was the most frequently observed among all deliveries (v2). However, this was the only one out of all 21 cases where the algorithm obtained results that were, on average, worse. With respect to loss, degradation cost was smaller by almost one-half for the algorithm than for manual assignment, while cost of wasted trucks was just one-third of this cost with manual assignment. However, for the cost of trucks left in the queue after the end of the day, was twice as large for the algorithm as for manual assignment. Nevertheless, this cost represented only 10\% of the total loss generated under the algorithm and 3\% under the manual method, thereby rendering its importance marginal.

With respect to preserving real-world grape-variety's distribution and intensity, the algorithm performed on average 3.85\% better than did the manual method. The highest result, 6.23\%, was obtained under the real Tuesday data frequency (fR), whereas the worst result, but still better than under the manual method, was generated for deliveries redistributed in frequent and moderate  four peaks (fP4). The highest loss for both decision methods was generated when deliveries were distributed into two high peaks (fP2), resulting from infrastructure overload generated by handling of many deliveries at one time. Similarly to variety-manipulation results, the algorithm obtained lower degradation and waste loss, while generating higher `left truck' loss.

The last tested aspect examined influence of a change in intensity, while preserving a realistic distribution of grape variety and a realistic frequency `shape.' Here the algorithm performed on average 5.21\% better than the manual method. Surprisingly, the smallest advantage (1.07\%) was obtained when the total number of deliveries was halved (i0.5) comparing to the real Tuesday data. In turn, the highest advantage (8.33\%) was generated for the other extreme, which was 1.5 times scaled the total number of real deliveries (i1.5). Again, due to the infrastructure overload, the highest loss was generated in the scenario with number of deliveries scaled at 1.5 times total number of actual deliveries.

As we have seen, the algorithm in general performed very well, exceeding payoffs generated by human intervention. However, in a weather-dependent industry, planned deliveries can drastically change due to unexpected weather conditions. For linear programming algorithms, such a change would make finding a maximizing strategy impossible. However, the same does not apply to dynamic programming solutions. After checking the algorithm’s performance when expected deliveries matched incoming deliveries, we tested the algorithm in cases of inconsistent delivery distribution.  Specifically, we tested inconsistency in each aspect: variety's distribution, frequency, and intensity, yielding 258 cases of inconsistency in total. Additionally, we averaged each result over four different realizations of the queue, which yielded 1032 simulations. The results from application of the algorithm when expected deliveries were inconsistent were arranged with results generated manually in such a way that both methods used the same queue.

Despite the inconsistent expected distribution, the algorithm generated 6.89\% better results compared to the manual method. A table containing all results can be reached at \cite{link}.  On average, the algorithm generated worse results in only 6 out of all 258 cases of inconsistency. The most interesting results were those for intensity inconsistency, where the actual queue intensity was 1.5 higher than the real-world intensity but the expected intensity was either 1 or 0.5 times scaled real-world intensity. The results in this case were much higher (15.7\% and 12.47\%, respectively) than for the consistent distribution (8.33\%), implying that, even if the expected delivery distribution were misspecified or were subject to sudden change, the algorithm in its basic version should on average generate better results than the manual method.

\section{Conclusion}\label{chapt_7}
Within this paper, we presented an optimization  algorithm (with respect to gain) for the real-time decision process of assigning trucks to pressing tanks in a specific wine factory. Based on our initial research, we concluded that the topic of real-time scheduling in the wine industry was not deeply explored, and, in fact, to date the most popular strategy for solving wine-connected scheduling problems employed linear programming (see \cite{wine_supply_chain_overview,OR_methods_overview_fresh_fruits}). However, the problem we seek to address requires to undertake series of decisions at real-time and so a linear programming solution scheme is not appropriate. Moreover, linear programming lacks robustness, a serious deficiency in the context of our problem. Consequently, we propose a new model that includes the stochasticity mentioned above. We also created an algorithm based on dynamic programming employing the Bellman equation. To our knowledge, this is the first attempt to apply dynamic programming to solving delivery scheduling problems originating in the wine industry.

Tests on real-world data demonstrated that, on average, it generates a higher payoff than manually assigning presses, especially in scenarios closest resembling the context within which these data were drawn. Due to its stochasticity, the algorithm can handle even misspecified delivery distributions and still generate a higher payoff than would the manual method.

The model, due to the promising results obtained in the course of this study, could receive serious consideration as a support software for pressing tank managers in wineries. Moreover, it can constitute a base solution method for further software development.

\begin{appendices}
\section{Distribution-consistent simulation results}\label{apx:DA_sim_results}
\begin{figure}[H]
\includegraphics[height=0.56\textheight, angle=90]{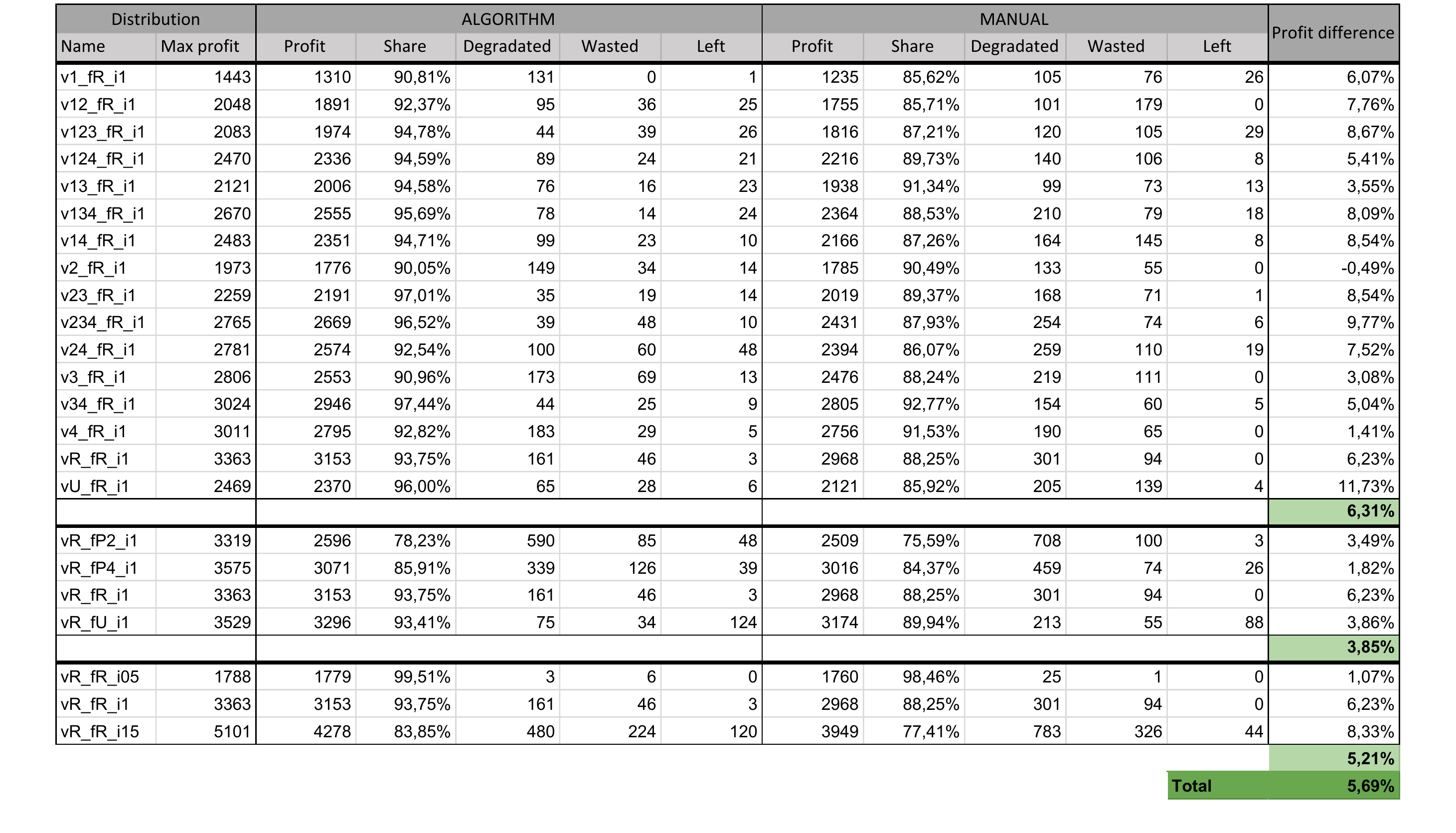}
\centering
\end{figure}

\section{Symbol list}
\renewcommand{\arraystretch}{1.6}
\begin{center}
\begin{longtable}{L{0.21\textwidth} L{0.6\textwidth} C{0.09\textwidth}}
\bfseries Sign & \bfseries Description & \bfseries Page\\
\hline
\endhead
\ensuremath{t} & Discrete unit of time, defined on the set $\{0,\ldots,T\}$. &  \hyperlink{page.3}{3} \\
\noalign{\global\arrayrulewidth=0.1mm}
  \arrayrulecolor{lightgray}\hline
\ensuremath{X_t} & State of the stochastic system at time $t$, i.e. description of the press at time $t$. &   \hyperlink{page.3}{3} \\
\ensuremath{\hspace{10pt}v(X_t)} & Variety of grapes in the press at time $t$. &  \hyperlink{page.4}{4} \\
\ensuremath{\hspace{10pt}l(X_t)} & Weight of grapes in the press at time $t$. &  \hyperlink{page.4}{4} \\
\ensuremath{\hspace{10pt}s(X_t)} & Processing start time of the press, checked at time $t$. &  \hyperlink{page.4}{4} \\
\ensuremath{\hspace{10pt}C} & Capacity of the press. &  \hyperlink{page.4}{4} \\
\ensuremath{\hspace{10pt}TP} & Total processing time of the press. &  \hyperlink{page.4}{4} \\
\noalign{\global\arrayrulewidth=0.1mm}
  \arrayrulecolor{lightgray}\hline
\ensuremath{Y_t} & Control at time $t$, i.e. decision which truck should be assigned to the press. &   \hyperlink{page.3}{3} \\
\ensuremath{\hspace{10pt}v(Y_t)} & Variety of grapes in the chosen truck $Y_t$. &  \hyperlink{page.4}{4} \\
\ensuremath{\hspace{10pt}l(Y_t)} &  Weight of grapes in the chosen truck $Y_t$. &  \hyperlink{page.4}{4} \\
\ensuremath{\hspace{10pt}a(Y_t)} & Arrival time of the chosen truck $Y_t$. &  \hyperlink{page.4}{4} \\
\ensuremath{\hspace{10pt}y_i} & The $i$th type of truck with the corresponding load variety and load weight $(v_i,l_i)$, $i=1,\ldots,N$. &  \hyperlink{page.4}{4} \\
\noalign{\global\arrayrulewidth=0.1mm}
  \arrayrulecolor{lightgray}\hline
\ensuremath{\Gamma(t,X_{t-1})} & Set of allowed controls on the state $X_{t-1}$ at time $t$, i.e. set of allowed truck types $y_i$ for the press at time $t$.&   \hyperlink{page.3}{3} \\
\noalign{\global\arrayrulewidth=0.1mm}
  \arrayrulecolor{lightgray}\hline
\ensuremath{Z_t} & New information arriving at time $t$, i.e. new trucks arriving at the winery. &   \hyperlink{page.3}{3} \\
\ensuremath{\hspace{10pt}Z_t^i} & Information about $i$th truck type presence among deliveries at time $t$.  &  \hyperlink{page.4}{4} \\
\ensuremath{\hspace{10pt}B_t^i} & Binary random variable indicating if the $i$th truck type will arrive to the winery at time $t$. &   \hyperlink{page.4}{4} \\
\noalign{\global\arrayrulewidth=0.1mm}
  \arrayrulecolor{lightgray}\hline
\ensuremath{f(t,X_{t-1},Y_t,Z_t)} & Transition function from the state $X_{t-1}$ to the state $X_t$, according to decision $Y_t$ and new information $Z_t$.  &   \hyperlink{page.3}{3} \\
\noalign{\global\arrayrulewidth=0.1mm}
  \arrayrulecolor{lightgray}\hline
\ensuremath{g(t,X_{t-1},Y_t)} & Payoff at time $t$ generated from the previous state $X_{t-1}$, according to the current decision $Y_t$.  &   \hyperlink{page.3}{3} \\
\noalign{\global\arrayrulewidth=0.1mm}
  \arrayrulecolor{lightgray}\hline
\ensuremath{V_{Y_{t\ldots T}}(t,x)} & Value function generated by $Y_{t\ldots T}$, i.e. expected sum of payoffs from time $t$  and state $x$, up to the time horizon $T$.  &   \hyperlink{page.3}{3} \\
\ensuremath{\hspace{10pt}V^*(t,x)} & Optimal value function generated by $Y^*_{t\ldots T}$, i.e. highest possible expected sum of payoffs.  &   \hyperlink{page.3}{3} \\
\ensuremath{\hspace{10pt}Y^*_{t\ldots T}} & Optimal strategy for the selection of controls. &  \hyperlink{page.3}{3} \\
\end{longtable}
\end{center}
\end{appendices}

\end{document}